\documentclass{amsart}
\usepackage{amssymb,amscd}

\newtheorem{lemma}{Lemma}
\newtheorem{prop}[lemma]{Proposition}
\newtheorem{cor}[lemma]{Corollary}
\newtheorem{thm}[lemma]{Theorem}

\newtheorem{thm?}[lemma]{Theorem?}

\newtheorem{conj}[lemma]{Conjecture}

\DeclareFontEncoding{OT2}{}{} 

\title{On the Hasse Principle for Shimura Curves}
\author{Pete L. Clark}

\address{1126 Burnside Hall \\ Department of Mathematics and Statistics \\
McGill University \\ 805 Sherbrooke West \\ Montreal, QC, Canada H3A 2K6}
\email{clark@math.mcgill.ca}

\newcommand{\F}{\ensuremath{\mathbb F}}

\newcommand{\Q}{\ensuremath{\mathbb Q}}
\newcommand{\R}{\ensuremath{\mathbb R}}
\newcommand{\Z}{\ensuremath{\mathbb Z}}

\newcommand{\C}{\ensuremath{\mathbb C}}

\newcommand{\ra}{\ensuremath{\rightarrow}}

\newcommand{\PP}{\ensuremath{\mathbb P}}

\newcommand{\Gal}{\operatorname{Gal}}

\newcommand{\End}{\operatorname{End}}

\newcommand{\lcm}{\operatorname{lcm}}

\newcommand{\OO}{\mathcal{O}}

\newcommand{\Jac}{\operatorname{Jac}}

\newcommand{\loc}{\operatorname{loc}}
\begin{document}
\maketitle

\begin{abstract}
Let $C$ be an algebraic curve defined over a number field $K$, of positive genus and without $K$-rational points.
We conjecture that there exists some extension field $L$ over which $C$ violates the Hasse principle, i.e., has points 
everywhere locally but not globally.  We show that our conjecture holds for all but finitely many Shimura curves of the 
form $X^D_0(N)_{/\Q}$ or $X^D_1(N)_{/\Q}$, where $D > 1$ and $N$ are coprime squarefree positive integers.  The proof uses a 
variation on a theorem of Frey, a gonality bound of Abramovich, and an analysis of local points of small degree.
\end{abstract}

\noindent
\section{Introduction}
\noindent
Fix $D > 1$ a squarefree positive integer, and $N \geq 1$ a
squarefree positive integer which is coprime to $D$.  Then there
exist Shimura curves $X^D_0(N)_{/\Q}$ (resp. $X^D_1(N)_{/\Q}$),
which are moduli spaces for abelian surfaces endowed with
endomorphisms by a maximal order in the indefinite rational
quaternion algebra $B_D$ of discriminant $D$ and $\Gamma_0(N)$
(resp. $\Gamma_1(N)$) level structure.  There are natural (``forgetful'')
modular maps $X^D_1(N) \ra X^D_0(N) \ra X^D = X^D_0(1)$.
\\ \indent
Understanding the points on these curves rational over
various number fields $K$ and their completions is an important
problem, analogous to the corresponding problem for the classical
modular curves $X_1(N), \ X_0(N)$ but even more interesting in at
least one respect. Namely, although the classical modular curves have
$\Q$-rational cusps, we have $X^D(\R) = \emptyset$. This raises the
possibility that certain Shimura curves, when considered over 
suitable non-real number fields $K$, may violate the
Hasse principle, i.e., may have points rational over every
completion of $K$ but not over $K$ itself.  And indeed, Jordan
showed that the curve $X^{39}_{/\Q(\sqrt{-23})}$ violates 
the Hasse principle \cite{Jordan} (see also \cite{Skoro}).  Later,
Skorobogatov and Yafaev provided explicit conditions on $D$, $N$
and $K = \Q(\sqrt{m})$ sufficient for $X^D_0(N)_{/K}$ to violate
the Hasse principle \cite{SY}.  Although it is plausible
that their conditions can be met for infinitely many choices of
$(D,N,K)$, these conditions (like Jordan's) include hypotheses
about the class group of $K$, so proving that they hold infinitely often 
seems very difficult.
\\ \\
Using different methods, we shall show that almost all 
Shimura curves of squarefree level violate the Hasse principle.  More precisely:
\begin{thm}
\label{MT1}
If $D > 546$, there is an integer $m$ such that $X^D_{/\Q(\sqrt{m})}$ 
violates the Hasse principle.
\end{thm}
\noindent

\begin{thm}
\label{MT2} There exists a constant $C$ such that if $D
\cdot N
> C$, then there exist number fields $K = K(D,N)$ and $L = L(D,N)$ such that
$X^D_0(N)_{/K}$ and $X^D_1(N)_{/L}$ violate the Hasse principle.
\end{thm}
\noindent
\begin{thm}
\label{MT3} Maintain the notation of the previous
theorem, and assume $D \cdot
N > C$. \\
a) We may choose $K$ such that $[K:\Q] \ | \ 4$. \\
b) Let $\{N_i\}$ be a sequence of squarefree positive integers
tending to infinity, and for each $i$, choose \emph{any}
squarefree positive integer $D_i > 1$ which is prime to $N_i$ and
such that $D_i \cdot N_i > C$. For all $i$, choose any number
field $L_i$ such that $X^{D_i}_1(N_i)_{/L_i}$ violates the Hasse
Principle.  Then $\lim_{i \ra \infty} [L_i:\Q] = \infty$.
\end{thm}
\noindent 
Remark: Using the same methods, one can show that there exists a positive 
integer $G$ such that if $F$ is a totally real field and $X_{/F}$ is a 
semistable Shimura curve of squarefree level of genus greater than $G$, then 
there exists an extension $K/F$ of degree dividing $4$ such that $X_{/K}$ 
violates the Hasse principle.  (Presumably the case of $\Gamma_1(N)$-level 
structure generalizes as well; I have not checked the details.)  
\\ \\
We choose to view these results as special cases of a more general
conjecture on algebraic curves defined over number fields.  Let
$V$ be a nonsingular, geometrically irreducible variety defined
over a number field $K$.  If there exists a number
field $L/K$ such that $V_{/L}$ has points rational over every
completion of $L$ but no $L$-rational points, we say that $V_{/K}$
is \textbf{a potential Hasse principle violation} (or
for brevity, ``$V_{/K}$ is PHPV.'')
\\ \\
Obviously $V_{/K}$ can only be PHPV if it has no $K$-rational
points.  Moreover, restricting to the class of curves, the case of
genus zero must be excluded, by Hasse-Minkowksi.  No further
restrictions spring to mind, and we conjecture that there are
none:
\begin{conj}
\label{MainConjecture}
Let $C_{/K}$ be a curve
defined over a number field $K$.  Assume that: \\
$(i)$ $C$ has no $K$-rational points. \\
$(ii)$ $C$ has positive genus. \\
Then $C_{/K}$ is a potential Hasse principle violation.
\end{conj}
\noindent A proof of this conjecture in the general case may not
be within current reach. However, using work of Faltings and Frey
we will derive a criterion for $C_{/K}$ to be PHPV (Theorem
\ref{Criterion}).  We use this criterion, together with an
analysis of local points on Shimura curves and a result of
Abramovich, to prove Theorems 1-3.
%
\\ \\
Finally, a warning: our method \textbf{does not} give an effective
procedure for finding the fields $K$ and $L$.  This is to be contrasted with the
work of \cite{Jordan} and \cite{SY}.
\\ \\
Remark and acknowledgment: I showed that there exist infinitely many 
Shimura curves $X^D_0(N)$ violating the Hasse principle over suitable quadratic fields in my 2003 Harvard thesis 
\cite[Main Theorem 5]{THESIS}.  Recently I became encouraged by the interest in this unpublished work shown by
V. Rotger, A. Skorobogatov and A. Yafaev \cite{RSY} and decided at last to write up the result.  Several improvements 
were found along the way, and the present results are stronger (and in some ways simpler) than what 
appeared in \cite{THESIS}.  I would like to take this opportunity -- better late than never -- to gratefully acknowledge 
the support of Harvard University and of my thesis advisor, Barry Mazur.
\section{Criteria for a Curve to be PHPV}
\noindent Let $X_{/K}$ be a variety over a field $K$. Define the
$m$-invariant $m(X) = m(X_{/K})$ to be the minimum degree of a
finite field extension $L/K$ such that $X(L) \neq \emptyset$.
\\ \\
If $K$ is a number field and $v$ is a place of $K$, we put $m_v(X)
:= m(X_{/K_v})$ and $m_{\loc}(X) = \lcm_v m_v(X)$.  Note that
applying Bertini's theorem and the Weil bound for curves over
finite fields, one gets the (well known) fact that $m_v(X) = 1$
for all but finitely many $v$,
so $m_{\loc}$ is well-defined.
\\ \\
The \emph{$K$-gonality} of a curve $C_{/K}$, denoted $d = d_K(C)$,
is the least positive integer $n$ for which there exists a degree
$n$ morphism $\varphi_{/K}: C \ra \PP^1$.  Clearly $d_K(C) = 1
\iff C \cong_K \PP^1$.  We say $C_{/K}$ is \emph{hyperelliptic} if
$d_K(C) \leq 2$.  (Beware that this is not quite the standard definition: under our 
definition all curves of genus zero and some curves of genus one are hyperelliptic, 
and a curve of higher genus may be hyperelliptic over $\overline{K}$ but not over $K$.)
\\ \\
For any curve $C_{/K}$ we have $m(C) \leq d_K(C)$.  Indeed, the
preimages of $\PP^1(K)$ under a degree $d = d_K(C)$ morphism
$\varphi: C \ra \PP^1$ yield infinitely many points $P$ on $C$ of
degree at most $d$.  Lying very much deeper is
the following result, which is a sort of converse.
\begin{thm}
\label{Frey} Let $C_{/K}$ be a curve over a number field, and, for
$n \in \Z^+$, let $\mathcal{S}_n(C)$ be the set of points $P \in
C(\overline{K})$ of degree dividing $n$.  If $\mathcal{S}_n(C)$ is
infinite, then $d_K(C) \leq 2n$.
\end{thm}
\noindent Proof: This is a small variation on a result of Frey
\cite[Prop. 2]{Frey} (which itself uses Faltings' spectacular theorem on rational points on 
subvarieties of abelian varieties).  In the original version, instead of
$\mathcal{S}_n(C)$ there appears the set $C^{(d)}(K)$ of points of
degree \emph{less than or equal to} $d$.  However, there also
appears the extra hypothesis that there exists $P_0 \in C(K)$,
which must be removed for our applications.  The existence of
$P_0$ is used to define a map from the $d$-fold symmetric product
$C^{(d)}$ to the Jacobian $\Jac(C)$, namely
\[\Phi: P_1 + \ldots + P_d \mapsto [P_1 + \ldots + P_d - dP_0]. \]
However, if $\mathcal{S}_n(C)$ is infinite, it is certainly
nonempty, so that for some $m \ | \ n$ there exists an effective
$K$-rational divisor $D_m$ of degree $m$, hence indeed an
effective $K$-rational divisor of degree $n$, namely $D_n =
\left(\frac{n}{m}\right) D_m$.  Then one can define the map
\[\Phi_D: P_1 + \ldots + P_d \mapsto [P_1 + \ldots + P_d -D_n], \]
and Frey's argument goes through verbatim with $\Phi_D$ in
place of $\Phi$.
%
\begin{thm}
\label{Criterion}
Let $C_{/K}$ be an algebraic curve defined over
a number field. Suppose: \\
a) $C(K) = \emptyset$. \\
b) $d_K(C) > 2m > 2$ for some multiple $m$ of $m_{\loc}(C)$. \\
Then there exist infinitely many extensions $L/K$ with $[L:K] =
m$ such that $C_{/L}$ is a counterexample to the Hasse
principle.
\end{thm}
\noindent Proof: By Theorem
\ref{Frey}, $\mathcal{S}_m(C)$ is a finite set. It follows that
the field $M$ defined as the compositum of $K(P)$ as $P$ ranges
through elements of $\mathcal{S}_m(C)$, is a number field. Let
$L/K$ be a number field which is linearly disjoint from $M/K$ and
such that $[L:K] = m$. Since $C(K) = \emptyset$, it follows that $C(L) = \emptyset$.
\\ \indent
Let $S = \{v_1, \ldots, v_r\}$ be the places of $K$ for which
$m_{v_i}(C) > 1$. By definition of $m_{\loc}(C) = \lcm_v m_v(C)$,
for each finite place $v_i \in S$, there exists a field extension
$L_i/K_{v_i}$ of degree ($m_{\loc}$, and \emph{a fortiori} of degree) 
$m$ such that $C(L_i)$ is nonempty.  At
each Archimedean place $v_i \in S$ (if any), we take $L_i$ to be
the $\R$-algebra $\C^{\frac{m}{2}}$.  Let $v_0$ be any finite
place not in $\Sigma$ and unramified in $M$, and let $L_0/K_{v_0}$
be a totally ramified extension of degree $m$.  For $0 \leq i \leq
r$ let $f_i \in K_{v_i}[x]$ be a defining polynomial for
$L_i/K_{v_i}$.  By weak approximation, for any $\epsilon> 0$,
there exists a degree $m$ polynomial $f \in K[x]$ such that, for
all $i$, each coefficient of $f-f_i$ has $v_i$-adic norm at most
$\epsilon$, so by Krasner's Lemma, for sufficiently small
$\epsilon$, $L = K[x]/(f)$ defines a degree $m$ field extension
with $L \otimes_K K_{v_i} \cong L_i$ for $0 \leq i \leq r$.  By
construction $m_{\loc}(C_{/L}) = 1$; moreover, $L/K$ is disjoint
from $M/K$, so $C(L) = \emptyset$.  By varying the choice of $v_0$
we clearly get infinitely many distinct fields $L$.
\\ \\
For $1 \leq n \leq 3$ there is a complete classification of
algebraic curves which have infinitely many points of degree at
most $n$ (Faltings for $n=1$, \cite{HS} for $n=2$,
\cite{AH} for $n=3$).  The quadratic case leads 
to the following ``supplement'' to Theorem
\ref{Criterion}.
\begin{thm}
\label{QuadraticCriterion}
 Suppose $C_{/K}$ is a curve over a number
field with $m_{\loc}(C) = 2$, $d_K(C) > 2$, and $C$ does not admit
a degree two morphism $\varphi: C \ra E$, where $E_{/K}$ is an
elliptic curve of positive rank.  Then there exist infinitely many
quadratic field extensions $L/K$ such that $C_{/L}$ is a
counterexample to the Hasse principle.
\end{thm}
\noindent
Proof: By the main result of \cite{HS}, the hypotheses imply that
$\mathcal{S}_2(C)$ is finite, and the rest of the proof is the
same as that of Theorem \ref{Criterion}.
\section{Local points on Shimura curves}
\noindent
We recall our notation: $D$ is the discriminant of a nonsplit indefinite rational 
quaternion algebra $B_D$, and is thus a nontrivial product of an even number of primes; $N$ 
is a squarefree positive integer prime to $D$.
\\ \\
The goal of this section is to prove the following result.
\begin{thm}
\label{MLOC}
a) For all $D$, $m_{\loc}(X^D) = 2$. \\
b) For all $D$ and all $N$, 
$m_{\loc}(X^D_0(N))$ is either $2$ or $4$.
\end{thm}
\noindent
Remark: Of course, for part b) one would like to know which of the two alternatives 
obtains.  It is possible to give a precise answer for this; more exactly, there are 
several nice sufficient conditions for $m_{\loc}(X^D_0(N)) = 2$ -- this holds, e.g., for 
fixed $D$ and all sufficiently large prime numbers $N$ -- and if none of these sufficient 
conditions hold there is a straightforward finite computation (coming from the Eichler-Selberg 
trace formula) that will determine the answer for any given pair $(D,N)$.  We hope to return to this and  
other matters pertaining to quadratic points on Shimura curves in a future work.

\subsection{Integral structure}
By work of Morita and Drinfeld there exists an integral canonical
model for $X^D_0(N)$.  In more words, there is a $\Z$-scheme,
projective, flat, and of relative dimension one, and which is a
coarse moduli space for Drinfeld's extension of the moduli problem
defining $X^D_0(N)_{\Q}$ to the category of schemes over $\Z$.
The generic fiber is canonically isomorphic to $X^D_0(N)$, so we
may denote this scheme by $X^D_0(N)_{\Z}$.  In particular, for
every prime number $p$ we have a special fiber $X^D_0(N)_{\F_p}$,
whose structure is highly relevant in the computation of
$m_p(X^D_0(N))$.  Depending upon $p$, the special fiber exhibits
three different kinds of behavior.
\subsection{Case of good reduction}
If $(p,DN) = 1$, $X^D_0(N)_{/\Z_p}$ is smooth.
\\ \\
As above, the special fiber $X^D_0(N)_{/\F_p}$ is again a moduli
space of abelian surfaces $A_{/k}$ (where $k$ is a field of characteristic $p$) endowed with a
quaternionic structure, i.e., an injection $\OO \hookrightarrow
\End A$. However, in contrast to characteristic zero, where the
generic QM surface is geometrically simple, \emph{all} QM surfaces
in characteristic $p$ are isogenous to $E \times E$, where
$E_{/k}$ is an elliptic curve \cite{Milne}.  Thus the full endomorphism
ring of $A$ depends on $\End(E)$: it is most often an order in
an imaginary quadratic field -- in which case we say $A$ (and its
corresponding point on the moduli space) is \emph{ordinary} -- but
there is a finite non-empty set of points for which the
endomorphism ring of $E$ is an order in a definite rational
quaternion algebra (of discriminant $p$).  The union of such
points determines the \emph{supersingular} locus on $X^D_{/\F_p}$,
and the supersingular locus on $X^D_0(N)_{/\F_p}$ is its complete
preimage under the modular ``forgetful'' map $X^D_0(N) \ra X^D$.
\begin{prop}
\label{Ogg}
There exists at least supersingular point on $X^D_0(N)_{\F_p}$.
Moreover, every supersingular point is defined over $\F_{p^2}$.
\end{prop} \noindent Proof: This is well known; see e.g. \cite{Ogg}.
\begin{cor}
\label{M1} For $p$ prime to $ND$, $m_p(X^D_0(N)) \leq 2$.
\end{cor}
\noindent Proof: By Proposition \ref{Ogg},
$m(X^D_0(N)_{\F_p}) \leq 2$.  Since $X^D_0(N)$ is smooth over
$\Z_p$, Hensel's Lemma implies that the $m$-invariant of the
generic fiber is at most the $m$-invariant of the special fiber,
and the result follows.
%
\subsection{Case of Cerednik-Drinfeld reduction}
For a positive integer $a$, we shall write $\Q_{p^a}$ for the unramified extension of $\Q_p$ of 
degree $a$ and $\Z_{p^a}$ for its valuation ring.
\\ \\
If $p \ | \ D$, then there exists a Mumford curve $C(D,N)_{/\Z_p}$
whose basechange to $\Z_{p^2}$ is isomorphic to
$X^D_0(N)_{/\Z_{p^2}}$.
\\ \\
Recall that a Mumford curve $C_{/\Q_p}$ is an algebraic
curve which can be uniformized by the $p$-adic upper half-plane.
Equivalently, the special fiber of the minimal regular model 
$\mathcal{C}_{/\Z_p}$ is a semistable curve which $\F_p$-split and degenerate: 
every irreducible component is isomorphic to $\PP^1$, with transverse 
intersections at $\F_p$-rational points. The special fiber is thus completely
determined by its \emph{dual graph}, a finite graph in which the
vertices correspond to the irreducible components $C_i$ and
the edges correspond to intersection points of $C_i$ and $C_j$.
Moreover, the degree of each vertex is at most $p+1$.
\begin{prop}
\label{M2}
 If $p \ | \ D$, $m_p(X^D_0(N)) \leq 2$.
\end{prop}
\noindent Proof: By Hensel's Lemma, it is enough to find a smooth
$\F_{p^2}$-rational point on $X^D_0(N) \cong_{\Z_{p^2}} C(D,N)$.
But any irreducible component $C_i$ on $C(D,N)_{\F_p}$ has at most
$p+1 = \#\PP^1(\F_p)$ singular points, hence at least $p^2+1-(p+1) > 0$ 
smooth $\F_{p^2}$-rational points.  The result follows.
\subsection{Case of Deligne-Rapoport reduction}
Suppose $p \ | \ N$.

\begin{thm} The arithmetic surface
$X^D_0(N)_{/ \Z_p}$ has the following structure: \\
a) The special fiber $X^D_0(N)_{\F_p}$ has two irreducible
components, each isomorphic to the (smooth) curve
$X^D_0(\frac{N}{p})_{\F_p}$.  \\
b) The two irreducible components intersect transversely at the
supersingular points, a corresponding point on the first copy of
$X^D_0(\frac{N}{p})_{\F_p}$ being glued to its image under the
Frobenius map.  \\
c) At each supersingular point $z \in X^D_0(N)_{/\F_{p^2}}$, the complete 
local ring is isomorphic to 
$\Z_{p^2}[[X,Y]]/(XY-p^{a_z})$ for some positive integer $a_z$.
\end{thm}
\noindent
Proof: The analogous statement for $D = 1$ (i.e., classical modular curves) is 
due to Deligne and Rapoport \cite{DR}.  In that same work it was pointed out that 
their result would continue to hold in the quaternionic context.  Apparently the first careful treatment of the 
quaternionic case is due to Buzzard \cite{Buzzard}, who worked however under the assumption of some additional 
rigidifying level structure.  The case of $X^D_0(N)$ (which is not a fine moduli scheme) was worked out 
(independently) in the theses of the author (unpublished) and of David Helm \cite[Appendix]{Helm}.


\begin{prop}
\label{M3} For any prime $p \ | \ N$, $m_p(X^D_0(N)) \leq 4$.
\end{prop}
\noindent Proof: Consider $X^D_0(N)_{\F_{p^2}}$.  Let $z$ be any
point on this curve coming from a supersingular point of
$X^D_0(\frac{N}{p})(\F_{p^2})$.  The local ring at $z$ is
analytically isomorphic to $\Z_{p^2}[[x,y]]/(xy-p^a)$ for some
integer $a \geq 1$.  \\ \indent Suppose first that $a > 1$.  Then
in order to get the minimal regular model, one must blow up the
point $z$ $(a-1)$ times, getting a chain of $a-1$ rational curves
defined over $\F_{p^2}$. Each of these curves has $p^2+1 -2$
smooth $\F_{p^2}$-rational points, which lift to give points on
$X^D_0(N)$ rational over $\Q_{p^2}$, the unramified quadratic
extension.
\\ \indent
Now assume that we have $a = 1$.  Let $K = \Q_{p^2}(\sqrt{p})$,
with valuation ring $R$ and uniformizing element $\pi = \sqrt{p}$.
 Then the completed local ring at $z$ of the integral model $X^D_0(N)_{/R}$
is isomorphic to $R[[x,y]]/(xy-\pi^{2})$.  We see that this local ring 
is no longer regular and must be blown up once to give a rational curve on the regular
model.  Thus, as in the previous case, we get a $K$-rational point
on $X^D_0(N)$, completing the proof.
\\ \\
Proof of Theorem \ref{MLOC}: From Propositions \ref{M1}, \ref{M2} and
\ref{M3}, we get $m_{\loc}(X^D_0(N)) \ | \ 4$.  Recalling that
$X^D_0(N)(\R) = \emptyset$, the result follows.

\section{Lower bounds on the gonality and genus}
\subsection{A result of Abramovich}
Let $\OO^1$ be the
group of norm $1$ units in a maximal $\Z$-order $\OO$ of $B_D$.
By restricting any embedding $B \hookrightarrow B \otimes_{\Q} \R
\cong M_2(\R)$ to $\OO^1$, one gets a realization of $\OO^1$ as a
discrete cocompact subgroup of $SL_2(\R)$, well-determined up to
conjugacy.  Then $\OO^1 \backslash \mathcal{H}$ is a compact
Riemann surface isomorphic to $X^D(\C)$.  For any finite index
subgroup $\Gamma \subset \OO^1$, we get a covering $X(\Gamma) \ra
X^D$.  
\begin{thm}(Abramovich)
\label{Abramovich}
\begin{equation} \frac{21}{200}(g(X^D(\Gamma))-1) \leq d_{\C}(X^D(\Gamma)).
\end{equation} 
\end{thm}

\subsection{Genus estimates}
In this section, all asymptotics are as $\min(D,N) \ra \infty$.  
\\ \\
For coprime squarefree positive integers $A$ and $N$, define
\[e_2(D,N) = \prod_{p \ | \ D} (1-(-4/p)) \prod_{q \ | \ N}
(1+(-4/q)),
\]
\[e_3(D,N) = \prod_{p \ | \ D} (1-(-3/p)) \prod_{q \ | \ N}
(1+(-3/q)).
\]
Let $\varphi$ be Euler's function, i.e., the multiplicative function such that 
$\varphi(p^k) = p^k - p^{k-1}$.  Let 
$\psi$ be the multiplicative function such that $\psi(p^k) = p^k+p^{k-1}$.
\\ \\
\noindent Then we have \cite[p. 301]{Ogg}:
\begin{equation}
\label{gDN} g(X^D_0(N)) = 1 +
\frac{1}{12}\varphi(D)\psi(N) - \frac{e_2(D,N)}{4}
- \frac{e_3(D,N)}{3} \sim \frac{\varphi(D)\psi(N)}{12}. 
\end{equation}
\noindent
In the case of $X^D_1(N)$ we will content ourselves with a lower bound for the genus.  Indeed, we will work 
instead with the Shimura covering $X^D_2(N)$, which is by definition the largest intermediate covering 
\[X^D_1(N) \ra Y \ra X^D_0(N) \]
such that $Y \ra X^D_0(N)$ is unramified abelian, say with Galois group $\Sigma(D,N)$.\footnote{We introduce 
the Shimura covering only to save a couple of lines of messy calculation with the Riemann-Hurwitz formula.}  From \cite[Corollary 1]{Ling}, we get 
\begin{equation}
\#\Sigma(D,N) \geq \frac{\varphi(N)}{2 \cdot \prod_{q | N} 6}. 
\end{equation}
Because $X^D_2(N) \ra X^D_0(N)$ is unramified, we have
\begin{equation}
g(X^D_2(N)) -1 =  \#\Sigma(D,N)(g(X^D_0(N))-1),
\end{equation}
and combining equations (2) through (4), we get
\[g(X^D_1(N)) -1 \geq g(X^D_2(N)) \gg \varphi(N) \cdot \frac{\varphi(D)}{24} \cdot \prod_{q \ | \ N} \frac{q+1}{6}. \]

\subsection{Proofs of the main theorems}
It is now only a matter of combining the results of $\S 2$, $\S 3$ and $\S 4$.
\\ \\
Consider first the case of $D = 1$.  According to Theorem \ref{MLOC}, $m_{\loc}(X^D) = 2$.  Moreover, by \cite{Rotger}, $X^D$ has only 
finitely many quadratic points when $D > 546$.  Thus Theorem \ref{QuadraticCriterion} applies to show 
that for such $D$, there exist infinitely many quadratic fields $K$ such that $X^D_{/K}$ violates the Hasse 
principle.
\\ \\
Next consider the case of $X^D_0(N)$.  By Theorem \ref{MLOC} we have $2 \ | \ m_{\loc}(X^D_0(N)) \ | \ 4$; 
in particular, $X^D_0(N)(\Q) = \emptyset$.  Applying Theorem \ref{Criterion} with $m=4$, we get that whenever the gonality 
$d_{\Q}(X^D_0(N)) > 8$, there exists infinitely many quartic fields $K/\Q$ such that $X^D_0(N)_{/K}$ violates 
the Hasse principle.  But by (1) and (2), $d_{\Q}(X^D_0(N)) \geq d_{\C}(X^D_0(N)) \gg \varphi(D) \psi(N)$.  
Since the latter quantity approaches infinity with $\min(D,N)$, the gonality condition holds for all but 
finitely many pairs $(D,N)$.
\\ \\
To deduce the $X^D_1(N)$ case, we need the following (essentially trivial) bound on the local m-invariant of 
$X^D_1(N)$.
\begin{lemma}
\label{EasyLemma}
Let $K$ be a number field and $f_{/K}: X_1 \ra X_2$ be a degree $d$ Galois covering of 
curves.  Then $m_{\loc}(X_1) \ | \  d \cdot m_{\loc}(X_2)$.
\end{lemma}
\noindent
Proof: First suppose that $f: X_2 \ra X_1$ is a Galois covering of curves defined over any field $K$, and let 
$P \in X_1(L)$ with $[L:K] = m(X_1)$.  Recall that the transitivity of the action of $\Gal(L(X_2)/L(X_1))$
on the points $\{Q_1,\ldots,Q_g\}$ of $X_2$ lying over $P$ implies the relation $efg = d$, where $e$ is the 
relative ramification index of $Q_i$ over $P$ and $f = [L(Q_i):L]$.  Thus $[L(Q_i):K] = f m(X_1) \ | \ d 
m(X_1)$.  The result follows easily by applying this observation at every place of $K$.
\begin{prop}
\[m_{\loc}(X^D_1(N)) \ | \ 2 \varphi(N). \]
\end{prop}
\noindent
Proof: The natural map $X^D_1(N) \ra X^D_0(N)$ is a Galois covering with group 
$(\Z/N\Z)^{\times}/(\pm 1)$.  Thus the result follows from Theorem \ref{MLOC} and Lemma \ref{EasyLemma}.
\\ \\
Now
\[d_{\Q}(X^D_1(N)) \geq d_{\C}(X^D_1(N)) \geq \frac{21}{200}(g(X^D_1(N)-1) \geq \frac{21}{200}(g(X^D_2(N))-1) \] \[ \gg
\varphi(N) \frac{\varphi(D)}{24} \prod_{q \ | N} \frac{q+1}{6}, \]
so that with $m = 2\varphi(N)$, we have
\[\frac{d_{\Q}(X^D_1(N))}{m} \gg \varphi(D) \prod_{q \ | \ N}\frac{q+1}{6}. \]
Again the right-hand side approaches infinity with $\min(D,N)$, so except for finitely many pairs $(D,N)$, 
we can find a field $L$ of degree $2\varphi(N)$ such that $X^D_1(N)_{/L}$ violates the Hasse principle.
\\ \\
Finally, Theorem 3b) is an immediate consequence of the following:
\begin{thm}
\label{CX}
For each prime number $p$ and positive integer $d \geq 1$, there exists a constant 
$N_0 = N_0(p,d)$ with the following property: for any $p$-adic field $K/\Q_p$ with 
$[K:\Q_p] \leq d$ and any integer $N \geq N_0$, $X^D_1(N)(K) = \emptyset$.
\end{thm}
\noindent
Proof: This is \cite[Theorem1]{CX}.
\\ \\
Indeed, taking, e.g., $p = 2$, we have that $m_2(X^D_1(N)) \ra \infty$ with $N$, uniformly in $D$.

\section{Final remarks}
\noindent
The constant $C$ in Theorem 1 can certainly be made explicit.  As the reader has probably noticed, 
it is easy to do so assuming a fixed number of prime divisors of $D \cdot N$, but slightly tedious in the 
general case.
\\ \\
More interesting than whittling down the set of excluded pairs $(D,N)$ to the optimal list which violate 
the inequality $d_{\Q}(X^D_0(N)) > 2m_{\loc}(X^D_0(N))$ is to investigate the cases of Conjecture 3 among 
Shimura curves of low (but positive!) genus.  For instance, suppose $X^D_0(N)$ has genus one.  
All such curves can be given in the form $y^2 = P(x)$, where $P \in \Q[x]$ has degree $4$ (and distinct 
roots).  For example,  $X^{14}$ is given by the equation $(x^2-13)^2 + 7^3 + 2y^2 = 0.$  It is not hard to see that 
$X^{14}_{/\Q(\sqrt{m})}$ has points everywhere locally if and only if $m$ is negative and prime to $7$, so that the set of such $m$ has 
density $\frac{3}{7}$ (as a subset of the set of all squarefree integers).  What can be said about the set of $m$ for which $X^{14}(\Q(\sqrt{m}))$ is nonempty?  
Or even about its density?  
\\ \\
One needs to know relatively little about the Shimura curves $X^D_0(N)$ beyond their semistability  to show that almost all of them 
satisfy Conjecture \ref{MainConjecture}.  Indeed, one can show the following
\begin{thm}
Let $\{X_n\}_{n=1}^{\infty}$ be a sequence of curves over a number field $K$.  \\ Suppose: \\
a) Each $X_n$ has everywhere semistable reduction. \\
b) $\lim_{n \ra \infty} \frac{d_K(X_n)}{\log g(X_n)} = \infty$. \\
c) There exists a fixed positive integer $A$ such that for all places $v$ and all $n$, the Galois action on the 
irreducible components of the special fiber $(X_n)_{/k_v}$ of the minimal model trivializes over an extension of degree $A$. \\
Then $X_n$ is PHPV for all sufficiently large $n$.
\end{thm}
\noindent
The hypotheses of theorem are satisfied (with $A = 2$) for the family of all semistable Shimura curves of squarefree level over any given totally 
real field $F$.  Note that condition b) looks very mild: surely a ``general'' family of curves, semistable or otherwise, should have gonality 
on the order of $g(X_n)$.  
\\ \\
I will save the proof for a later time, at which point I hope it will be more clear whether condition c) can be weakened or removed.

\end{document}